\documentclass[12pt]{article}

\newtheorem{theo}{Theorem}[section]
\newtheorem{lem}[theo]{Lemma}

\newtheorem{defi}[theo]{Definition}

\newcommand{\coLim}{{\rm colim}}
\newcommand{\CC}{{\bf C}}
\newcommand{\DD}{{\,\bf D}}

\newcommand{\Z}{{\mathbb Z}}
\newcommand{\F}{{\mathcal F}}

\usepackage{amssymb,latexsym}
\def\leq{\leqslant}

\def\geq{\geqslant}

\title
{
Flows in Graphs and Homology of Free Categories
~\footnote
{Research is supported by RF grant center at Novosibirsk
State University in Russia and by T\"UB\.ITAK and NATO in Turkey}
}
\author{Ahmet A. Husainov, Hamza \c{C}al{\i}\c{s}{\i}c{\i}}
\date{}

\begin{document}

\maketitle

\begin{abstract}
We introduce the notion of a generalized flow 
on a graph with coefficients in a $R$-representation 
and show that the module 
of flows is isomorphic to the first derived functor of the colimit. 
We generalize Kirchhoff's laws and build an exact sequence 
for calculating the module of flows on the union of graphs.
\end{abstract}



\section{Preliminaries}

This work is devoted to the $R$-module of generalized flows in a graph.
We applicate the homology of small categories in the sense of \cite{hil1971}. 
Our approach is distinguished with the theory described 
in Wagner's work \cite{wag1998} where groups of flows are considered 
as a new invariant of graphs.

Let $\CC$ be a small category, $R$ a ring with identity. 
Denote by $Mod_R$ the category of left $R$-modules and $R$-homomorphisms,
$Mod_R^\CC$ the category of functors $\CC \rightarrow  Mod_R$, 
$\coLim^\CC:  Mod_R^\CC \rightarrow Mod_R$ the colimit functor. 
The category  $Mod_R^\CC$  has enough projectives \cite{hil1971}.
The functor $\coLim^\CC$ is right exact.
Hence for every integer $n\geq 0$ it is defined  
the $n$-th right derived functor
$\coLim_n^\CC:  Mod_R^\CC \rightarrow Mod_R$ of the colimit.
Let $F:\CC \rightarrow Mod_R$ be a functor. 
For arbitrary family $\{a_i\}_{i\in I}$ we will say that 
{\em almost all $~a_i$ are zeros} if there exists a finite 
subset $J\subseteq I$ such that $a_i= 0$  for all $i\in I\setminus{J}$.
Denote 
$C_n(\CC,F) = \sum\limits_{c_0\rightarrow ... \rightarrow c_n}F(c_0)$ 
and write elements of $C_n(\CC,F)$
as sums 
$\sum\limits_{c_0\rightarrow ... \rightarrow c_n}
f_{c_0\rightarrow ... \rightarrow c_n}[c_0\rightarrow ... \rightarrow c_n]$
with  $f_{c_0\rightarrow ... \rightarrow c_n} \in F(c_0)$ where almost all 
$f_{c_0\rightarrow ... \rightarrow c_n}$ are zeros.
For every $c_0\rightarrow ... \rightarrow c_{n+1}$ and $f\in F(c_0)$
we let 
$$
d_n(f [c_0\rightarrow ... \rightarrow c_{n+1}]) = 
F(c_0\rightarrow c_1)(f)
[c_1\rightarrow ... \rightarrow c_{n+1}] +
$$
$$
\sum\limits_{i=1}^{n+1}(-1)^i f[c_0\rightarrow ... \rightarrow \hat{c_i} 
\rightarrow ... \rightarrow c_{n+1}]
$$ 
and 
define homomorphisms  $d_n: C_{n+1}(\CC, F)\rightarrow C_n(\CC,F)$ by
$$
d_n\left(\sum\limits_{c_0\rightarrow ... \rightarrow c_{n+1}}
f_{c_0\rightarrow ... \rightarrow c_{n+1}}
[c_0\rightarrow ... \rightarrow c_{n+1}]\right) = 
$$
$$
\sum\limits_{c_0\rightarrow ... \rightarrow c_{n+1}}
d_n\left(f_{c_0\rightarrow ... \rightarrow c_{n+1}}
[c_0\rightarrow ... \rightarrow c_{n+1}]\right)
$$
Here
$$
[c_0 \stackrel{\alpha_1}{\rightarrow}
\cdots \stackrel{\alpha_i}{\rightarrow} \hat{c_i}
\stackrel{\alpha_{i+1}}\rightarrow \cdots
\stackrel{\alpha_{n+1}}\rightarrow c_{n+1}
] = 
$$

$$
\left\{
\begin{array}{ll}
[c_0 \stackrel{\alpha_1}\rightarrow 
\cdots \stackrel{\alpha_{i-1}}\rightarrow c_{i-1} 
\stackrel{\alpha_{i+1}\alpha_i}\rightarrow
c_{i+1}\stackrel{\alpha_{i+2}}\rightarrow \cdots
\stackrel{\alpha_{n+1}}\rightarrow c_{n+1}], & ~if~ 0<i<n+1,\\

[c_0 \stackrel{\alpha_1}\rightarrow 
\cdots 
\stackrel{\alpha_{n}}\rightarrow c_{n}], &  ~if~i=n+1.\\

\end{array}
\right.
$$

It is well-known \cite[Application 2]{gab1967} that $R$-modules $\coLim_n^\CC{F}$
are isomorphic to homologies of the complex
$$
0 \stackrel{d_{-1}}\longleftarrow C_0(\CC,F) \stackrel{d_0}\longleftarrow C_1(\CC,F) 
\stackrel{d_1}\longleftarrow \cdots \stackrel{d_{n-1}}\longleftarrow C_n(\CC,F)  
\stackrel{d_n}\longleftarrow\cdots
$$
in the sense of 
$\coLim_n^\CC F \cong Ker d_{n-1}/Im d_n$, $n\geq 0$. 

If  
$\alpha\circ\beta = id$ in $\CC$ implies $\alpha=id$ and $\beta=id$,
then $\CC$ is called to be {\em without retractions}.
If $\CC$  is a small category without retractions, 
then the complex $\{C_n(\CC,F), d_n\}$ includes the subcomplex 
$\{C_n^+(\CC,F), d_n\}$  where  
$$
C_n^+(\CC,F) = \sum\limits_{c_0{\rightarrow\atop\not=}
\cdots{\rightarrow\atop\not=}c_n} F(c_0)
$$ is the submodule in which 
sequences $c_0\rightarrow\cdots \rightarrow c_n$ do not contain 
identity morphisms if $n>0$, with $C_0^+(\CC,F)=C_0(\CC,F)$.
As the dual affirmation \cite[Proposition 2.2]{X1997} we can prove the following.

\begin{lem}\label{L11} 
Let $\CC$  be a small category without retractions, 
$R$ a ring with identity. Then for each functor 
$F: \CC \rightarrow Mod_R$  the $R$-modules 
$\coLim_n^\CC F$   are isomorphic to $n$-th homology modules 
of the complex $\{C_n^+(\CC,F), d_n\}$.
\end{lem}

Let $\Delta_\CC\Z$ be a functor from a small category to the category 
of Abelian groups and homomorphisms which assign to every 
$c \in \CC$ the group 
of integers $\Z$ and to every $\alpha \in Mor\CC$ the identity 
homomorphism  $id_\Z: \Z\rightarrow\Z$. We denote by  $H_n(\CC)$ the groups 
$\coLim_n^\CC\Delta_\CC\Z$ for all $n\geq 0$.
Let $S:\CC \rightarrow \DD$ be a functor between small categories, 
$d$ an object in $\DD$. The {\em comma-category} $~d/S$ is defined as 
the following category:

Objects of $d/S$ are pairs $(c,\alpha)$ with $c\in Ob\CC$ and 
$\alpha \in \DD(d,S(c))$, morphisms 
$(c_1,\alpha_1)\rightarrow (c_2,\alpha_2)$ in $d/S$ consist of the 
triples $(\beta\in\CC(c_1,c_2), \alpha_1, \alpha_2)$ satisfying 
$S(\beta)\circ\alpha_1 = \alpha_2$. A functor $S:\CC \rightarrow \DD$ 
is called {\em strong cofinal} if $H_n(d/S)\cong H_n(pt)$ for all $n\geq 0$. 
Here $pt=\{*\}$ is the discrete category with one object, thus 
$S$ is strong cofinal if and only if for all $d\in Ob\DD$ the groups $H_n(d/S)$ 
are zeros for all $n>0$ and  $d/S$ are connected.

By Oberst's Theorem \cite[Theorem 2.3]{obe1968} if $S$ is strong cofinal, 
then for every functor $F:\DD\rightarrow Mod_R$  the canonical 
homomorphisms $\coLim_n^\DD F \rightarrow \coLim_n^\CC (F\circ S)$ are 
isomorphisms.

Let $\CC$ be a small category. Each object $a\in Ob\CC$ will be 
considered  as the identity $id_a$, so $Ob\CC \subseteq Mor\CC$. 
The {\em factorization category} \cite{bau1985} ${\F}\CC$ is the 
category such that $Ob {\F}\CC = Mor\CC$ and for every pair 
$f, g \in Mor\CC$ the set ${\F}\CC(f,g)$ of 
morphisms consists of the pairs $(\alpha,\beta)$, with 
$\alpha,\beta \in Mor\CC$, 
for which $\beta\circ{f}\circ\alpha = g$.

\begin{lem} \label{L12}
The functor $s: ({\F}\CC)^{op}\rightarrow\CC$ which assign to any $f\in Ob {\F}\CC$ 
its domain $s(f)$ and to any morphism $(\alpha,\beta)$ the morphism $\alpha$
is strong cofinal.
\end{lem}
{\sc Proof.} Objects of $c/s$ for $c\in Ob\CC$ are pairs $(x,\alpha)$ 
of morphisms 
$c\stackrel{x}\rightarrow s(\alpha)\stackrel\alpha\rightarrow{t(\alpha)}$, 
and morphisms $(x,\alpha)\rightarrow(y,\beta)$ are commutative diagrams
$$
\begin{array}{ccccc}
c & \stackrel{x}{\longrightarrow} &
s(\alpha) & \stackrel{\alpha}{
\longrightarrow} & t(\alpha)\\
\downarrow\lefteqn{id_c} &&
\downarrow && \uparrow\\
c & \stackrel{y}{\longrightarrow}&
s(\beta) & \stackrel{\beta}{\longrightarrow} & t(\beta)
\end{array}
$$
For every $c\in Ob\CC$ the category $c/s$ includes the full subcategory 
consisting of all the objects $(id_c,\alpha)$. 
This subcategory is isomorphic to $(c/\CC)^{op}$. For every object 
$(x,\alpha)$ there is a morphism 
$(id_c,\alpha\circ{x})\rightarrow (x,\alpha)$ 
such that for each morphism  
$(id_c,\beta)\rightarrow (x,\alpha)$
there exists the unique morphism 
$(id_c,\beta)\rightarrow (id_c,\alpha\circ{x})$
for which the following diagram is commutative
$$
\begin{array}{ccc}
(id_c,\alpha\circ{x}) & {\longrightarrow} &
(x, \alpha) \\
\uparrow\lefteqn{\exists{!}} &&
 \uparrow\\
(id_c,\beta) &  \stackrel{id}{\longrightarrow} & (id_c,\beta)
\end{array}
$$
It follows that there exists a right adjoint functor to the inclusion 
$(c/\CC)^{op}\subseteq c/s$. A right adjoint is strong cofinal, hence 
$$
	\coLim_n^{c/s}\Delta\Z \cong \coLim_n^{(c/\CC)^{op}}\Delta\Z.
$$
But $(c/\CC)^{op}$ has a terminal object. Thus $H_n(c/s)\cong H_n(pt)$  
for all $n\geq 0$. $\Box$

\section{Generalized Flows}
	
By a {\em (directed) graph}  we mean a pair of sets $(A,V)$  
and a pair of functions 
$A{{{{s}\atop\rightarrow}}\atop{{\rightarrow\atop{t}}}} V$. 
The elements of $A$ are called {\em arrows}\,, $V$ is the set of 
{\em vertexes}\,,  $s(\alpha)$ and  $t(\alpha)$ are called the 
{\em source} and the {\em target} of $\alpha\in A$ respectively.

	Let  $\Gamma = (A,V,s,t)$ be a graph, $R$ a ring with identity. 
A {\em $R$-representation} of $\Gamma$ is a family of $R$-modules 
$\{F(v)\}_{v\in V}$ with a family of homomorphisms 
$\{F(\alpha): F(s(\alpha))\rightarrow F(t(\alpha))\}_{\alpha\in A}$.
A {\em path} in $\Gamma$ from $u\in V$ to $v\in V$ is an arbitrary 
word  
$\alpha_1\alpha_2\cdots\alpha_n$ with $t(\alpha_1)=v$ 
and $s(\alpha_n)=u$ such that $s(\alpha_i)=t(\alpha_{i+1})$ for 
all $1\leq i \leq n-1$. 
For each 
vertex $v\in V$ define  $id_v$ as the empty path from $v$ to $v$.
Objects of the {\em category of paths in $\Gamma$} are vertexes $v\in V$, 
and morphisms are paths in $\Gamma$ with the composition law 
$$
\alpha_1\cdots\alpha_n\circ\beta_1\cdots\beta_m = 
\alpha_1\cdots\alpha_n\beta_1\cdots\beta_m, 
$$
for $s(\alpha_n) = t(\beta_1)$.
Let $W\Gamma$  be the category of paths in $\Gamma$.
For every $R$-representation $F$ there is the unique 
functor $\tilde{F}: W\Gamma \rightarrow Mod_R$ such that 
$\tilde{F}\vert_{\Gamma} = F$.  It is defined by 
$\tilde{F}(\alpha_1\cdots\alpha_n) = F(\alpha_1)\cdots{F}(\alpha_n)$ if 
$n>0$, and $\tilde{F}(id_v) = id_{F(v)}$.

\begin{defi} Let $\Gamma$ be a graph, $F$ a 
$R$-representation of $\Gamma$. A {\em flow on $\Gamma$ 
with coefficients in $F$} is a family  $\{f_\gamma\}_{\gamma\in A}$ 
of elements $f_\gamma \in F(s(\gamma))$ such that 
almost all of $f_\gamma$ are zeros and 
for each $v\in V$ the following equality holds:
$$
\sum\limits_{s(\gamma)=v}f_\gamma = \sum\limits_{t(\gamma)=v}F(\gamma)(f_\gamma)
$$
\end{defi}
We denote by $\Phi(\Gamma; F)$ the $R$-module of flows on $\Gamma$ 
with coefficients in $F$. We have the following exact sequence:
\begin{equation}\label{seq01}
0 \rightarrow \Phi(\Gamma;F) \rightarrow \sum_{\gamma\in A}F(s(\gamma))
\stackrel{d}\longrightarrow \sum_{v\in V}F(v) \rightarrow 
\coLim^{W\Gamma}\tilde{F} \rightarrow 0 ,
\end{equation}
where $ d(\sum\limits_{\gamma\in A}f_\gamma [\gamma])_v =
\sum\limits_{t(\gamma)=v}F(\gamma)(f_\gamma) -
\sum\limits_{s(\gamma)=v}f_\gamma$.

Let $S:\CC \rightarrow \DD$ be 
a functor between small categories, $d\in Ob\DD$ an object. 
We denote by $S/d$ the category which objects are pairs $(c,\alpha)$
with
$c\in Ob\CC$, $\alpha\in\DD(S(c),d))$; morphisms 
$(c_1,\alpha_1)\rightarrow(c_2,\alpha_2)$ in $S/d$ are triples 
$(\beta\in\CC(c_1,c_2),\alpha_1,\alpha_2)$ satisfying 
$\alpha_2\circ{S(\beta)} = \alpha_1$. It is clear that 
$d/(S^{op}) \cong (S/d)^{op}$. Denote ${\F}W\Gamma = {\F}(W\Gamma)$.
\begin{lem}\label{L22}
Let $\Gamma=(A,V,s,t)$ be a graph, $\Gamma'$ the full subcategory 
of ${\F}W\Gamma$  such that $Ob\Gamma' = A \coprod V$. 
Then ${\Gamma'}^{op}$  is strong cofinal in $({\F}W\Gamma)^{op}$.
\end{lem}
{\sc Proof.} We denote by $S$ the inclusion $\Gamma' \subseteq {\F}W\Gamma$. 
The objects of $S/w$ are pairs $(x,(\alpha,\beta))$ 
where $\alpha$ and $\beta$ are paths in $\Gamma$ satisfying
$\beta\circ{x}\circ\alpha=w$ with either $x=id$ or $x\in A$. 
Hence for 
$w = (v_0\stackrel{\alpha_1}\rightarrow v_1\stackrel{\alpha_2}\rightarrow
\cdots \stackrel{\alpha_n}\rightarrow v_n)$ the category $S/w$ 
is the following:
\begin{center}
\begin{picture}(180,120)
\thicklines
\put(20,10){The~ category~$S/w$}

\put(0,45){\llap{$f_1$}}
\put(40,45){\llap{$f_2$}}
\put(0,60){\circle*{3}}
\put(40,60){\circle*{3}}
\put(120,60){\circle*{3}}
\put(120,100){\circle*{3}}
\put(160,100){\circle*{3}}
\put(120,45){$f_n$}
\put(70,80){\circle*{3}}
\put(80,80){\circle*{3}}
\put(90,80){\circle*{3}}
\put(0,110){\llap{$e_0$}}
\put(40,110){\llap{$e_1$}}
\put(0,100){\circle*{3}}
\put(40,100){\circle*{3}}
\put(160,110){$e_n$}

\put(0,100){\vector(0,-1){38}}
\put(40,100){\vector(-1,-1){38}}
\put(40,100){\vector(0,-1){38}}
\put(60,80){\vector(-1,-1){18}}
\put(120,100){\vector(0,-1){38}}
\put(160,100){\vector(-1,-1){38}}
\put(100,80){\line(1,1){20}}

\end{picture}
\end{center}
with the objects $e_0 =
(id_{v_0},(id_{v_0},\alpha_n\cdots\alpha_1))$,  $f_1= 
(\alpha_1,(id_{v_0},\alpha_n\cdots\alpha_2))$,  $e_1 = 
(id_{v_1}, (\alpha_1, \alpha_n\cdots\alpha_2))$, 
 $\cdots$, 
$f_n = (\alpha_n,(\alpha_{n-1}\cdots\alpha_1,id_{v_n}))$, 
$e_n = (id_{v_n}, (\alpha_n \cdots \alpha_1, id_{v_n}))$, 
where all morphisms excluding $e_{i-1}\rightarrow f_i$ and 
$e_i\rightarrow f_i$, for $i \in \{1, 2, \cdots , n\}$,  are identities.

But $w/S^{op}$ is isomorphic to $(S/w)^{op}$, hence  $H_n(w/S^{op})=0$
for $n >0$ and $w/S^{op}$ is connected. $\Box$

\begin{lem}\label{L23}
Let $\Gamma = (A, V, s, t)$ be a graph, 
$F: ({\F}W\Gamma)^{op} \rightarrow Mod_R$ a functor. 
Then $\coLim_1^{({\F}W\Gamma)^{op}}F$ is isomorphic to the submodule 
in $\sum\limits_{\gamma\in A}F(\gamma)$  consisting of families  
$\{g_\gamma\}_{\gamma\in A}$ for which the following equality holds 
for every $v\in V$:
\begin{equation}\label{eqvg}
\sum_{v=s(\gamma)}F(v\stackrel{(id,\gamma)}\rightarrow\gamma)g_\gamma = 
\sum_{v=t(\gamma)}F(v\stackrel{(\gamma,id)}\rightarrow\gamma)g_\gamma
\end{equation}
Here $\gamma$ runs the set $A$.
\end{lem}
{\sc Proof.} We have by Lemma \ref{L22} from the Oberst theorem 
that 
$\coLim_n^{({\F}W\Gamma)^{op}}F$ is isomorphic to 
$\coLim_n^{{\Gamma'}^{op}}F\vert_{{\Gamma'}^{op}}$.
The category ${\Gamma'}^{op}$ has no retractions except identities. 
Hence, by Lemma \ref{L11} $R$-modules 
$\coLim_n^{({\F}W\Gamma)^{op}}F$ are isomorphic to the homology groups 
of the chain complex 
$$
0\rightarrow C_1^+({\Gamma'}^{op},F)\rightarrow C_0^+({\Gamma'}^{op},F)
\rightarrow 0
$$
which is equal to
$$
	0 \rightarrow \sum_{id\not=(v\rightarrow\gamma)\in Mor\Gamma'}F(\gamma)
	\stackrel{d}\longrightarrow \sum_{\gamma\in A}F(\gamma) 
	\oplus \sum_{v\in V}F(v) \rightarrow 0
$$
with 
$$
d(\sum\limits_{v\rightarrow\gamma}f_{v\rightarrow\gamma}[v\rightarrow\gamma]) 
= 
\sum\limits_{v\rightarrow\gamma}d(f_{v\rightarrow\gamma}[v\rightarrow\gamma])
=
\sum\limits_{v\rightarrow\gamma}(F(v\rightarrow\gamma)(f_{v\rightarrow\gamma})
[v] - 
f_{v\rightarrow\gamma}[\gamma]).
$$
(Here we consider the homology of the category which is opposite to 
$\Gamma'$, in particular $f_{v\rightarrow\gamma} \in F(\gamma)$ and 
a homomorphism $F(v\rightarrow\gamma)$ acts from $F(\gamma)$ into $F(v)$.) 
Therefore $\coLim_1^{{\Gamma'}^{op}} F$ is isomorphic to the 
$R$-module of families
$f_{v\rightarrow\gamma}\in F(\gamma)$
satisfying $f_{s(\gamma)\rightarrow\gamma}+f_{t(\gamma)\rightarrow\gamma}=0$                                        
for each $\gamma\in A$, and 
$\sum\limits_{v\rightarrow\gamma}F(v\rightarrow\gamma)(f_{v\rightarrow\gamma})
= 0$ for each $v\in V$.
	
	We denote $g_\gamma = f_{s(\gamma)\rightarrow\gamma}$  
for a such family. Then $f_{t(\gamma)\rightarrow\gamma}= -g_\gamma$, and  
$\sum\limits_{v=s(\gamma)}F(v\rightarrow\gamma)(g_{\gamma}) =
\sum\limits_{v=t(\gamma)}F(v\rightarrow\gamma)(g_{\gamma})$. 
Thus $\coLim_1^{{\F}W\Gamma^{op}}F$ is isomorphic to the submodule 
of $\sum\limits_{\gamma\in A}F(\gamma)$  consisting of 
$\{g_\gamma\}_{\gamma\in A}$ for which the equation (\ref{eqvg}) holds.
$\Box$

\begin{theo}\label{T21} Let $\Gamma = (A,V,s,t)$   be a graph, $R$ a ring 
with identity, $F$ a $R$-representation of $\Gamma$. Then   
$\Phi(\Gamma;F) \cong \coLim_1^{W\Gamma}\tilde{F}$.
\end{theo}
{\sc Proof.} By Lemma \ref{L12} the functor $s: ({\F}\CC)^{op}\rightarrow\CC$ 
is strong cofinal for an arbitrary small category. Hence 
$\coLim_1^{W\Gamma}\tilde{F} \cong \coLim_1^{{\F}W\Gamma^{op}}(\tilde{F}\circ{s})$. 
The substitution of $\tilde{F}\circ{s}$ instead $F$ 
in  Lemma \ref{L23} leads to 
the concluding that $\coLim_1^{W\Gamma}\tilde{F}$ 
is isomorphic to the submodule of $\sum\limits_{\gamma\in A}F(\gamma)$   
which consists of  families  $f_\gamma\in F(s(\gamma))$ satisfying 
$\sum\limits_{v=s(\gamma)}f_\gamma = 
\sum\limits_{v=t(\gamma)}F(\gamma)(f_\gamma)$, 
for each $v\in V$.  $\Box$


\section{ The First Kirchhoff Law}

Let $\Gamma = (A, V, s, t)$   be a graph, $F$ a $R$-representation 
of $\Gamma$. Elements of $\sum\limits_{v\in V}F(v)$  and  
$\sum\limits_{\gamma\in A}F(s(\gamma))$  are called {\em $0$-chains} 
and {\em $1$-chains} respectively. 

Let  
$\varepsilon: \sum\limits_{v\in V}F(v) \rightarrow \coLim^{W\Gamma}\tilde{F}$ 
be the canonical $R$-homomorphism. It follows from  the exact 
sequence (\ref{seq01}) that the equation $d{f} = \varphi$ has a solution 
for $\varphi \in \sum\limits_{v\in V}F(v)$ if and only if 
$\varepsilon(\varphi) = 0$. Hence $\coLim^{W\Gamma}\tilde{F}$ can be 
interpreted as the $R$-module of "obstructions". Denote it by
$\Phi_0(\Gamma;F)$. We have the exact sequence
$$
0 \rightarrow \Phi(\Gamma;F) \stackrel{\subseteq}\rightarrow 
\sum_{\gamma\in A}F(s(\gamma))
\stackrel{d}\longrightarrow \sum_{v\in V}F(v) 
\stackrel{\varepsilon}\rightarrow 
\Phi_0(\Gamma; F) \rightarrow 0 
$$
with  $\Phi(\Gamma;F) \cong \coLim_1^{W\Gamma}\tilde{F}$ and
 $\Phi_0(\Gamma;F)\cong \coLim^{W\Gamma}\tilde{F}$.
      
A {\em network} $(\Gamma, E, F)$ consists of the following data:

1) a graph $\Gamma = (A, V, s, t)$;

2) an arbitrary subset $E \subseteq V$ 
which elements are called {\em external};

3) a $R$-representation $F$ of $\Gamma$.

We say that $1$-chain  $\{f_\gamma\}_{\gamma\in A}$  {\em satisfies to the 
first Kirchhoff law} if  
$$
	d (\{f_\gamma\}_{\gamma\in{A}})_v = 0, \quad \forall v \not\in E .
$$
Let $\Phi(\Gamma,E;F)$ be a $R$-module of all 
$1$-chains satisfying to the first Kirchhoff law in the network 
$(\Gamma, E, F)$.
For $E=\emptyset$ an $1$-chain satisfies to the first Kirchhoff law if 
and only if it is a flow on $\Gamma$ with coefficients in $F$.
Thus, $\Phi(\Gamma,\emptyset;F) = \Phi(\Gamma;F)$.

A vertex $v\in V$ is called to be {\em attractive} if there are not 
arrows with $s(\gamma)=v$.
Let $E\subseteq{V}$ be any subset such that all $e\in E$ are 
attractive vertexes,
$F_E$ the $R$-representation of $\Gamma$ with 
$F_E(v) = 0$ for all $v\not\in E$,
and $F_E(v)=F(v)$ for all $v\in E$.
We have by Theorem \ref{T21} the following

\begin{lem}
Let $(\Gamma, E, F)$ be a network. If all vertexes in $E$ are 
attractive then
$\Phi(\Gamma, E; F) \cong \coLim_1^{W\Gamma}(\widetilde{F/F_E})$.
\end{lem}

To a  description the $R$-module of $1$-chains satisfying to the  
first Kirchhoff law in any network $(\Gamma, E, F)$ we add to the 
graph $\Gamma$ the vertex $*$ and  the arrows $\gamma_e$ for all
$e\in E$ with $s(\gamma_e)=e$ and $t(\gamma_e)=*$. Denote by 
$\Gamma\cup_E{pt}$ the obtained graph. Let $F\oplus_E{0}$ be the 
$R$-representation of $\Gamma\cup_E{pt}$ such that 
$(F\oplus_E{0})\vert_\Gamma = F$ and $(F\oplus_E{0})(*) = 0$.

\begin{theo} 
For any network $(\Gamma, E, F)$ the $R$-module $\Phi(\Gamma, E; F)$ 
is isomorphic to $\coLim_1^{W(\Gamma\cup_E{pt})}\widetilde{(F\oplus_E{0})}$.
\end{theo}
{\sc Proof.} Consider the network $(\Gamma\cup_E{pt},{pt},F\oplus_E{0})$. 
The vertex $*$ is attractive in $\Gamma\cup_E{pt}$. It follows from the  
previous lemma that 
$\Phi(\Gamma\cup_E{pt},{pt};F\oplus_E{0}) \cong 
\coLim_1^{W(\Gamma\cup_E{pt})}\widetilde{(F\oplus_E{0})}$. But 
$\Phi(\Gamma\cup_E{pt},{pt};F\oplus_E{0}) = \Phi(\Gamma, E; F)$. $\Box$


\section{Flows on the Union of Graphs}

Let $(I,\leq)$ be  a partially ordered set. A covering 
$X = \bigcup\limits_{i\in I}X_i$ of a set  is called to be 
{\em locally filtered}
if $i<j$ in $I$ implies $X_i \subseteq X_j$   and if for each  
$x \in X_i \cap X_j$ there exists $k\in I$ such that $k < i$, $k<j$, 
and $x \in X_k \subseteq X_i \cap X_j$.

	A graph $\Gamma = (A, V, s, t)$ is called to be 
{\em locally filtered covered} by graphs $\Gamma_i = (A_i, V_i, s_i, t_i)$  
if 
$A=\bigcup\limits_{i\in I}A_i$ and $V=\bigcup\limits_{i\in I}V_i$ 
are locally filtered coverings and
the following diagrams are commutative 
$$
\begin{array}{ccc}
~A_i & {\subseteq}\atop
{\longrightarrow} & A_j\\
{s_i}\lefteqn\downarrow && {s_j}\lefteqn\downarrow\\
~V_i & {\subseteq}\atop
{\longrightarrow} & V_j
\end{array}
\quad
\begin{array}{ccc}
~A_i & {\subseteq}\atop
{\longrightarrow} & A_j\\
\downarrow\lefteqn{t_i} && \downarrow\lefteqn{t_j}\\
~V_i & {\subseteq}\atop
{\longrightarrow} & V_j
\end{array}
$$
for all $i\leq j$ in $I$. 
\begin{theo}
Let $\Gamma=(A,V,s,t)$ be a graph which is locally filtered 
covered by graphs   
$\{\Gamma_i\}_{i\in I} = \{A_i, V_i, s_i, t_i\}$, $E \subseteq V$ 
a subset such that $E = \bigcup\limits_{i\in I}E_i$ is a 
locally filtered covering by $E_i\subseteq V_i$. 
Then for each $R$-representation $F$ of 
$\Gamma$ there exists an exact sequence
$$
0 \rightarrow \coLim_2^I\{\Phi_0(\Gamma_i,E_i; F_i)\} 
	\rightarrow \coLim^I\{\Phi(\Gamma_i,E_i; F_i)\}
$$
$$
	 \rightarrow \Phi(\Gamma,E; F) \rightarrow 
	    \coLim_1^I\{\Phi_0(\Gamma_i,E_i; F_i)\} \rightarrow 0
$$
where $F_i = F\vert_{\Gamma_i}$.
\end{theo}
{\sc Proof.} At first we consider the case $E=\emptyset$. 
Let  $\Gamma'_i\subseteq \Gamma'$ are the categories defined for 
the graphs $\Gamma_i\subseteq \Gamma$ in Lemma \ref{L22}. 
The covering  $\{\Gamma'_i\}_{i\in I}$  satisfies 
to the conditions of \cite[Corollary 3.2]{X1989}. 
The functors $s:({\F}W\Gamma)^{op}\rightarrow W\Gamma$ and 
$S^{op}:{\Gamma'}^{op}\rightarrow ({\F}W\Gamma)^{op}$ are strong cofinal by
Lemma \ref{L12} and Lemma \ref{L22}. Hence, $\coLim_1^{W\Gamma}\tilde{F}$ 
is isomorphic to 
$\coLim_1^{{\Gamma'}^{op}}(\tilde{F}\circ{s}\circ{S^{op}})$.
By \cite[Corollary 3.2]{X1989} there is the spectral sequence with 
$E_{p,q}^2 = \coLim_p^I\{\coLim_q^{{\Gamma'}_i^{op}}\tilde{F}\circ{s}
\circ{S^{op}}\vert_{{\Gamma'}_i^{op}}\}$ 
which converges to 
$\coLim_{n}^{{\Gamma'}^{op}}\tilde{F}\circ{s}\circ{S^{op}}$. 
The substitution $\Gamma_i$ instead $\Gamma$ leads to the functors
$s_i:({\F}W\Gamma_i)^{op}\rightarrow W\Gamma_i$ and 
$S_i^{op}:{\Gamma'_i}^{op}\rightarrow ({\F}W\Gamma_i)^{op}$ which 
are strong cofinal by Lemmas \ref{L12} and  \ref{L22}. 
It follows from $(\tilde{F}\circ{s}\circ{S^{op}})\vert_{{\Gamma_i'}^{op}}
= \tilde{F}\vert_{W\Gamma_i}\circ{s_i}\circ{S_i^{op}}$ and from the 
strong confinality of $s_i$ and $S_i^{op}$ that we have the 
spectral sequence
$\coLim_p^I\{\coLim_q^{W\Gamma_i}\tilde{F}\vert_{W\Gamma_i}\}
\Rightarrow \coLim_{p+q}^{W\Gamma}\tilde{F}$.
Then the exact sequence of terms of low degree \cite[P.332]{mac1975} 
gives the exact sequence
$$
0 \rightarrow \coLim_2^I\{\Phi_0(\Gamma_i; F_i)\} 
	\rightarrow \coLim^I\{\Phi(\Gamma_i; F_i)\}
$$
$$
	 \rightarrow \Phi(\Gamma; F) \rightarrow 
	    \coLim_1^I\{\Phi_0(\Gamma_i; F_i)\} \rightarrow 0
$$

For $E\not=\emptyset$ we consider the locally filtered covering 
$\Gamma\cup_E pt = \bigcup\limits_{i\in I}(\Gamma_i\cup_{E_i}pt)$.
There is an exact sequence
$$
0 \rightarrow \coLim_2^I\{\Phi_0(\Gamma_i\cup_{E_i}pt; F_i\oplus_{E_i}0)\} 
	\rightarrow \coLim^I\{\Phi(\Gamma_i\cup_{E_i}pt; F_i\oplus_{E_i}0)\}
$$
$$
	 \rightarrow \Phi(\Gamma\cup_E{pt}; F\oplus_E{0}) \rightarrow 
	    \coLim_1^I\{\Phi_0(\Gamma_i\cup_{E_i}pt; F_i\oplus_{E_i}0)\} 
	    \rightarrow 0
$$
The equalities 
$\Phi_0(\Gamma\cup_E{pt}; F\oplus_E{0}) = \Phi_0(\Gamma, E; F)$ and 
$\Phi(\Gamma\cup_E{pt}; F\oplus_E{0}) = \Phi(\Gamma, E; F)$
give 
looking. $\Box$

\section{ The Second Kirchhoff Law}

Let $R$ be a field.
The {\em internal product} on a vector space  $T$ is a bilinear 
map $<,>: T \times T \rightarrow R$ such that $<a,b> = <b,a>$ for 
all $a, b \in T$. 
A network $(\Gamma, E, F)$ together with an internal 
product $<,>$ on $\sum\limits_{\gamma\in A}F(s(\gamma))$
is called to be {\em Euclidian} if the implication 
$<f,f>=0 \Rightarrow f=0$ is true. 
We say that an $1$-chain 
$f = \{f_\gamma\}$ of an Euclidian network satisfies to the 
{\em second Kirchhoff law} if the linear map $<f,->: 
\sum\limits_{\gamma\in A}F(s(\gamma)) \rightarrow R$ has zero values 
on $\Phi(\Gamma,F)$ in the sense that 
$$
<f,->\vert_{\Phi(\Gamma,F)} = 0.
$$
\begin{theo}
Let $(\Gamma,E,F)$ be an Euclidian network in which $R$ is a field,
$\Gamma$ a finite graph, and $F(v)$ finite dimensional vector spaces for 
all $v\in V$. Then for each $\varphi\in \sum\limits_{v\in E}F(v)$ 
satisfying $\varepsilon(\varphi) = 0$ there is the unique $1$-chain 
$f = \{f_\gamma\}$ such that $df = \varphi$ and 
$<f,->\vert_{\Phi(\Gamma,F)} = 0$.
\end{theo}
{\sc Proof.} If $df=0$ then $f\in \Phi(\Gamma,F)$, in this case 
$<f,->\vert_{\Phi(\Gamma,F)}=0$ implies $<f,f>=0$ and $f=0$. 
Considering $f_i$, with $i\in\{1,2\}$, for which $df_i = \varphi$ and 
$<f_i,->\vert_{\Phi(\Gamma,F)}=0$, we obtain 
$f_1-f_2 = 0$. Hence, the solution is unique.

Consider a map 
$\sum\limits_{\gamma\in A}F(s(\gamma))\stackrel{\eta}\rightarrow 
Mod_R(\Phi(\Gamma,F),R) \oplus \varepsilon^{-1}(0)$ which assign to any 
$g\in \sum\limits_{\gamma\in A}F(s(\gamma))$ the pair 
$<g,->\vert_{\Phi(\Gamma,F)}\oplus dg$. It is easy to verify that 
this map is the isomorphism. We let $g =\eta^{-1}(0\oplus\varphi)$.
Then $\eta(g)= 0\oplus\varphi$. Hence $dg = \varphi$ and 
$<g,->\vert_{\Phi(\Gamma,F)} = 0$. Thus there exists the solution. 
$\Box$

\begin{flushleft}
Department of Computer Technologies,\\
Komsomolsk-on-Amur State Technical University,\\
prosp. Lenina, 27, Komsomolsk-on-Amur, Russia, 681013\\
and\\
Amasya Egitim Fakultesi,\\
Matematik Bolumu, Ondokuz Mayis University,\\
Amasya, Turkey, 05189\\
Email: husainov@knastu.ru and hmcalisici@yahoo.com
\end{flushleft}

\end{document}